\documentclass[12pt]{article}
\input{psfig}
\usepackage{amssymb}

 \newtheorem{theorem}{Theorem}[section] 
 
 \newtheorem{proposition}[theorem]{Proposition}

 \def\Box
  {\hfill \thinspace\vbox{\hrule height .5pt \hbox{\vrule  
   width .5pt \vbox to 7pt{\hbox to 3.5pt{}} \vrule width .5pt} 
   \hrule height 0pt depth .5pt}}

 \newenvironment{proof}{{\it Proof:\/}}{$\Box$\vskip 0.08in}
 
\title{Topological Insights from the Chinese Rings}

\author{J\'ozef H. Przytycki and Adam S. Sikora}

\date{}
\begin{document}

\maketitle

\begin{abstract}
L. Kauffman conjectured that a particular solution of the Chinese Rings
puzzle is the simplest possible. We prove his conjecture by using 
low-dimensional topology and group theory. We notice also a
surprising connection between the Chinese Rings and Habiro moves
(related to Vassiliev invariants).
\end{abstract}


\section{Introduction}
In some of the popular puzzles one is supposed to take a ring off a rope 
which is usually tangled with the rigid part of the puzzle. Typically, 
such puzzles possess ingenious solutions which, if considered carefully 
enough, lead to interesting problems in low-dimensional topology.
This phenomenon can be observed in the Chinese Rings
which are shown below. The purpose of the puzzle is to take the loose 
ring off the rope. \vspace*{.2in}

\centerline{\psfig{figure=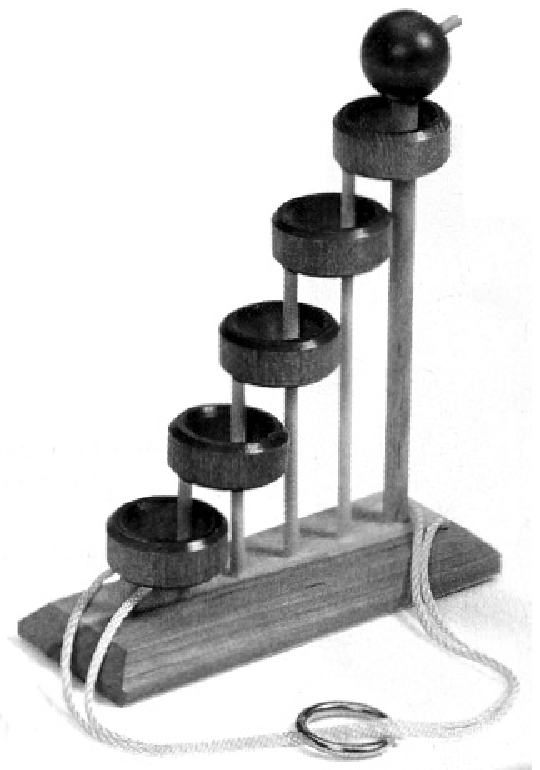,height=5cm}}

Although this puzzle was for a long time a toy of the first author's
children, we become seriously interested in it only after L. Kauffman 
pointed to us some interesting mathematical aspects of it, \cite{K2}. 
The history of the Chinese Rings can be found in \cite{K2,BC}.
The goal of this paper is to prove that a particular solution
of the Chinese Rings puzzle is the simplest possible, as conjectured in 
\cite{K2}.

If you have not thought about this puzzle before, then we suggest that
you try to solve it now in order to gain an appreciation for this
beautiful problem.
Do not be discouraged by the initial difficulty in finding a solution.
Indeed it is almost impossible to see all the necessary moves to be 
performed on the rope just by staring at the picture above.
However, it is easy to see that a solution exists: Imagine for a moment
that the toy is made out of a flexible material and press the
longer columns down so that all columns have the same length, 
as shown in Figure 1.
\vspace*{.2in}

\centerline{\psfig{figure=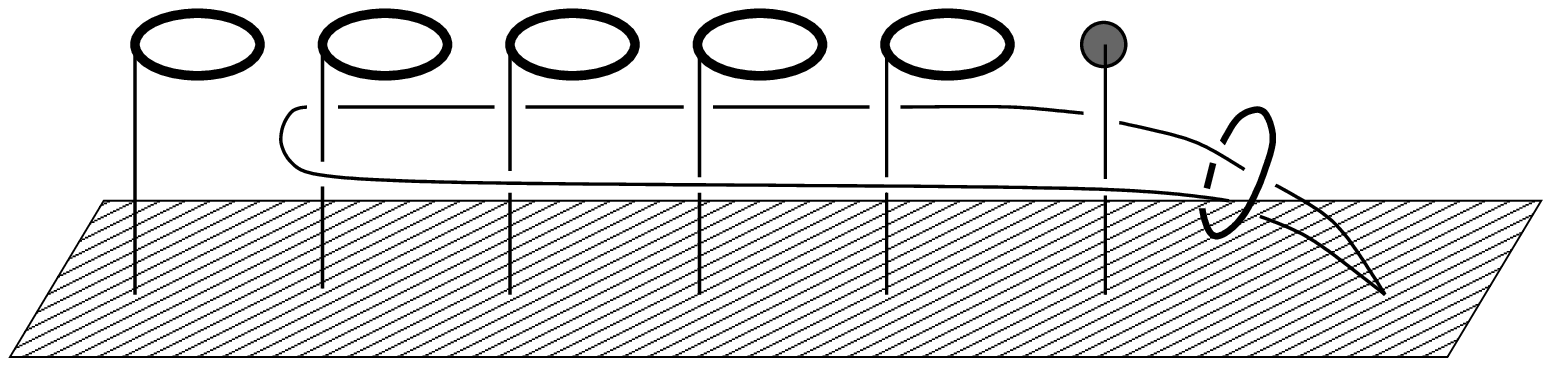,height=2.5cm}}
\begin{center}
                Figure 1. 
\end{center}

Now the solution becomes obvious!
Since now the rope is completely separated from the columns,
it can also be untangled in the original
puzzle. (Actually, we need to be sure that
the rope is long enough for this theoretical solution to be correct; 
this can be confirmed by direct experimentation).

The above solution shows how a simple topological idea can yield
a beautiful solution to what seems to be a complicated problem.
The same idea shows that the Chinese Rings can be presented in a somewhat
more regular form:

\centerline{\psfig{figure=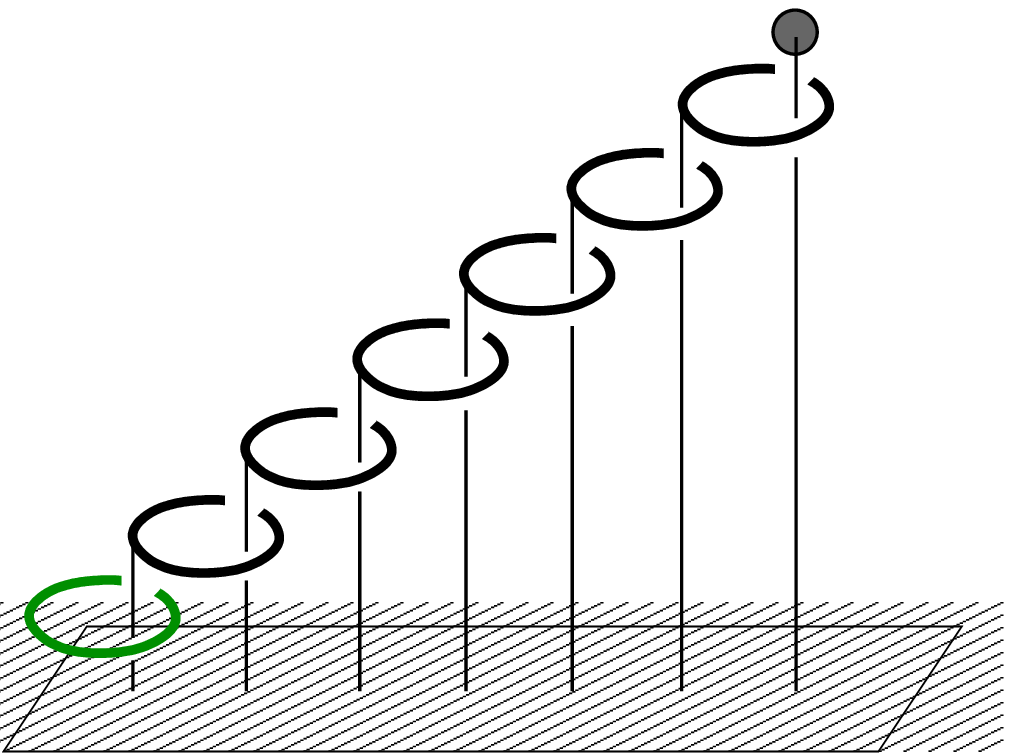,height=3cm}}
\begin{center}
                Figure 2. 
\end{center}

Note that we replaced the rope by another column with an attached ring.
The problem now is to separate the loose ring (which is assumed to be 
infinitely elastic) from the solid part of the puzzle.
The problem addressed in this paper is to find the simplest solution to
this puzzle.
As we will see soon, the solution to this problem involves an inductive
argument in the number of columns. For that reason we consider a
more general version of the puzzle of Fig. 2 composed of an arbitrary
number of columns (all but one with an attached ring at the top).

A precise answer to our problem requires an objective measure of 
the complexity of possible solutions to the puzzle. For that, imagine 
an arc $A,$ drawn below with a dashed line,
joining the highest column of the toy with the base. 
The complexity of a solution for the puzzle is the
minimal number of times the elastic ring passes through the arc $A$ in the
process of that solution.

\centerline{\psfig{figure=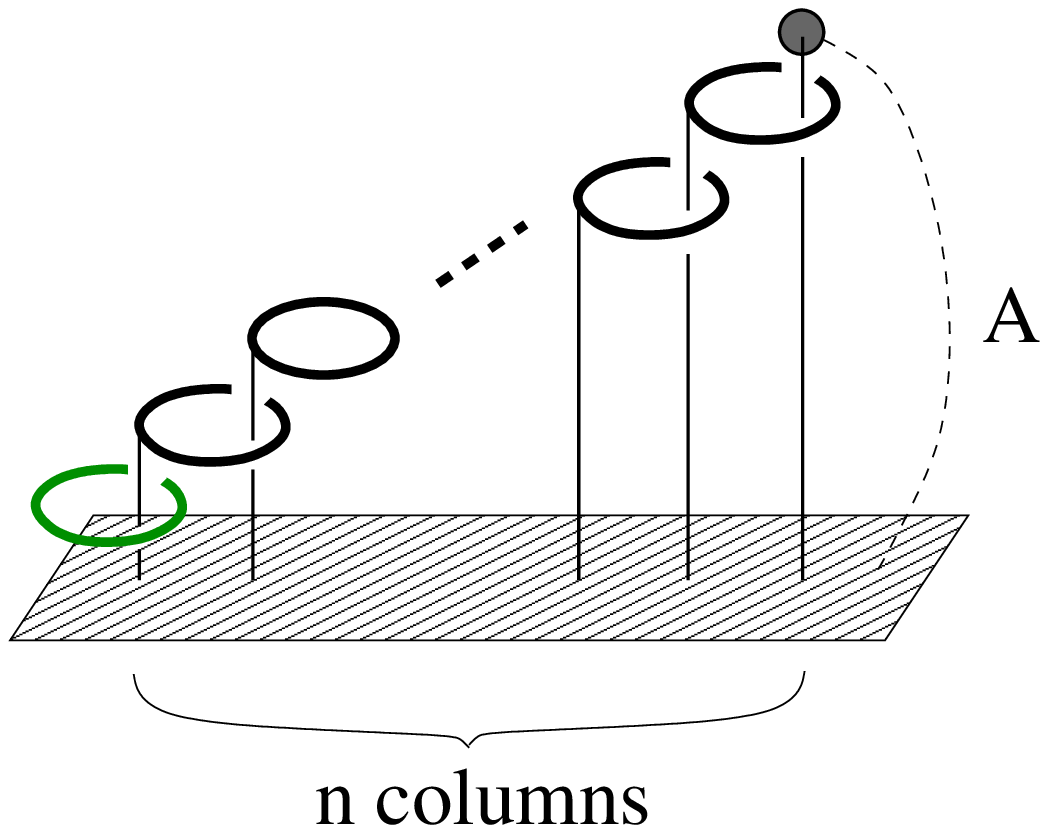,height=3cm}}
\begin{center}
                Figure 3. 
\end{center}


\begin{theorem}
The minimal complexity of a solution for the Chinese Rings with
$n$ columns presented at Fig 3 is $2^{n-1}.$ 
\end{theorem}

This result was conjectured by L. Kauffman in \cite{K2}, who wrote:

"This problem in the topology of the Chinese Rings is a useful test case for
questions that can arise in applications of knot theory to natural structures
where there is always a mixture of topology and mechanical/geometrical
modeling.

A solution to the Ring Conjecture will probably involve the discovery
of new techniques for understanding topology of graph embeddings in
three-dimensional space. It is fun to be able to take a classical
puzzle as fascinating as the Chinese Rings  and find within it
a significant topological problem. Let us find the solution!
"

{\bf Habiro Moves}\ There is a surprising connection between the
Chinese Rings (known for hundreds of years) and a very recent theory of
Vassiliev invariants of knots, \cite{BL}.
K. Habiro \cite{Ha} proved that two knots cannot be distinguished by 
Vassiliev invariants of order $\leq n$ if and only if they are related
by a sequence of moves (called Habiro moves) presented below:\\
\ \\
\ \\
\centerline{\psfig{figure=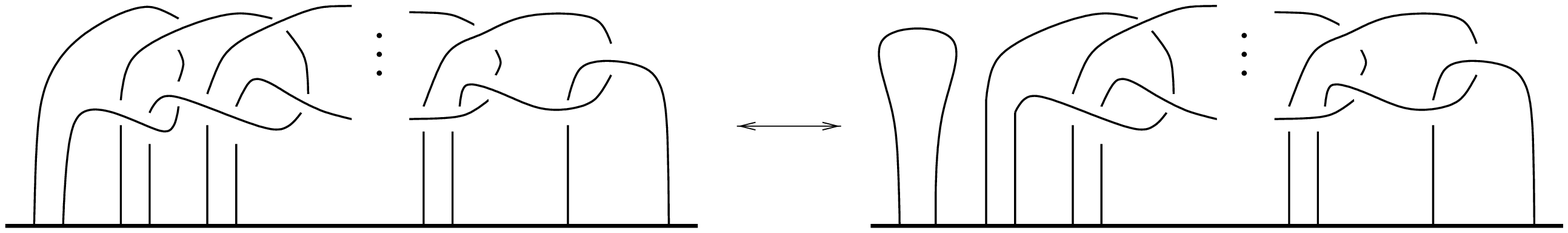,height=2cm}}
\ \\

Note that the Habiro move corresponds exactly to the operation of
untangling the loose ring in the version of Chinese Rings presented
in Figure 3!

\section{Proof}

The solution of the puzzle presented at the beginning of the paper 
requires that we first untangle the rings attached to the columns.
Note that as a result of this deformation of the puzzle, the arc $A$
will assume the form presented below.\vspace*{.2in}

\centerline{\psfig{figure=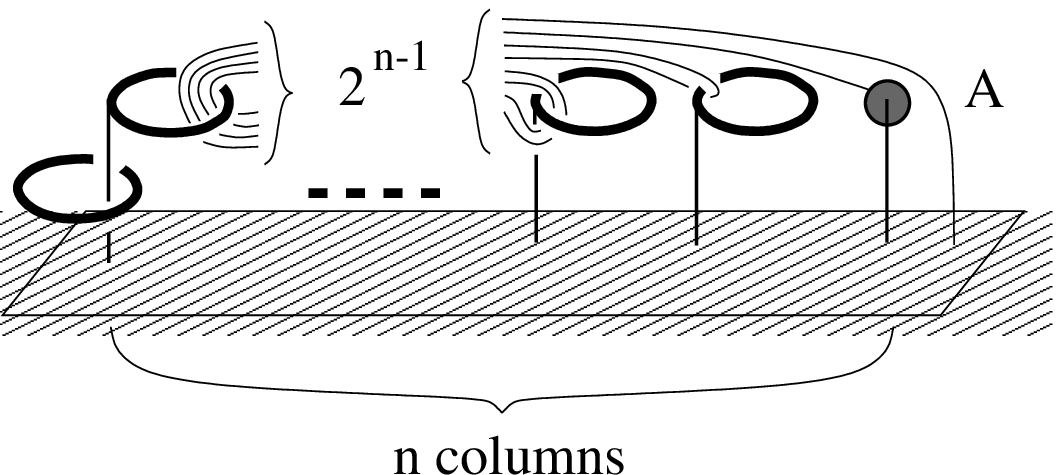,height=3.2cm}}
\begin{center}
                Figure 4. 
\end{center}

Because we can now remove the ring by passing through the $2^{n-1}$
strands of $A,$ the complexity of this algorithm is obviously at most
$2^{n-1}.$
We are going to show that this is the simplest solution to the
Chinese Rings; i.e. the complexity of any other solution is
not less than $2^{n-1}.$ For the purpose of the proof we are allowed
to relax the conditions of the problem and consider it in a topological
setting by allowing arbitrary deformations of the dimensions
and the shape of the toy. Moreover we add a point at infinity to the 
ambient three-space and hence consider the problem in $S^3.$ 
These assumptions make the problem simpler and
surely do not increase the complexity of the minimal solution.

The body of the Chinese Rings with $n$ columns (and with the
loose ring excluded) is a handlebody $H_{n-1}$ of genus $n-1$
embedded in the standard way into $S^3.$ Its complement,
$H_{n-1}',$ is also a handlebody of 
genus $n-1.$ The loose ring is contractible in $H_{n-1}',$ and this is 
the reason for which the puzzle has a solution. The arc $A$ assumes 
a complicated position in $H_{n-1}'$ which can be deduced from 
Fig. 4, see Fig. 5\vspace*{.2in}

\centerline{\psfig{figure=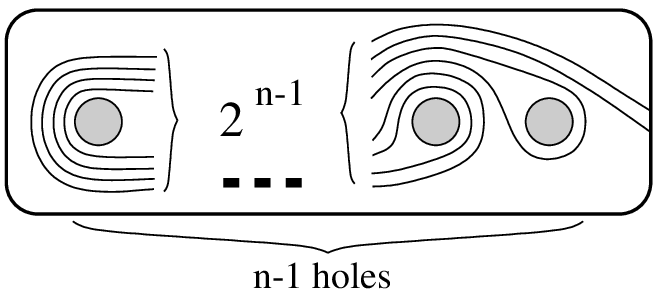,height=2.7cm}}
\begin{center}
                Figure 5: Complement of the puzzle $=H_{n-1}'$  
\end{center}

If we remove the arc $A$ from $H_{n-1}'$ then the loose ring, $S^1,$
will no longer be contractible. Recall that the complexity of a
solution for the puzzle is the minimal number of passes of the ring
through $A$ necessary for
contracting the ring to a point in $H_{n-1}'\setminus A.$ 
We need to show that this number is at least $2^{n-1}.$
Unfortunately, the position of the ring $S^1$ in $H_{n-1}'\setminus A$ is very
complicated; try to figure it out by yourself to see that, indeed, it is 
not an easy problem! In order to avoid this problem we use three tricks:
First, we take a dual approach:
we fix the ring $S^1$ in $H_{n-1}'$ and count the number of times
the arc $A$ has to pass through $S^1$ in order to make $S^1$ contractible
in $H_{n-1}'\setminus A.$
Second, we slightly relax the conditions of the puzzle, by assuming that
$A$ can be deformed by an arbitrary homotopy fixing its endpoints.
We assume that the endpoints of $A$ are at some $x_0\in \partial H_{n-1}'.$
Therefore we consider $A$ as an element of $\pi_1(H_{n-1}'\setminus S^1,x_0).$
Our intention is to prove that under these relaxed conditions the
complexity of the puzzle is at least $2^{n-1}.$ This will surely imply that
the complexity of the puzzle under original assumptions is also at least
$2^{n-1}.$

Third, we modify the puzzle by
attaching the loose ring to the base of the puzzle by an additional column, 
and we try to untangle $A$ in the modified puzzle.
This seems to be a more difficult problem, but actually it is not.
Observe that the modified puzzle is homeomorphic to a handlebody of
genus $n,\ H_n$ and that its complement is homeomorphic also to a handlebody,
$H_n'.$ Notice also that we have an embedding 
$i: H_n'\to H_{n-1}'\setminus S^1$ corresponding to the fact
that the modified toy differs from the original only by an additional
column. One can prove using Van-Kampen's Theorem that the map $i$
induces an isomorphism of the fundamental groups, $i_*:
\pi_1(H_n')\to \pi_1(H_{n-1}'\setminus S^1).$ Therefore the additional 
column does not create any new obstacle for $A$! (Recall that we
consider the arc $A$ up to homotopy only.) The benefit of this trick is that
the modified puzzle is much easier to solve because it does not have
any loose ring. Observe that the arc $A$ lies in $H_n'$ in the pattern
presented in Figure 6.

\centerline{\psfig{figure=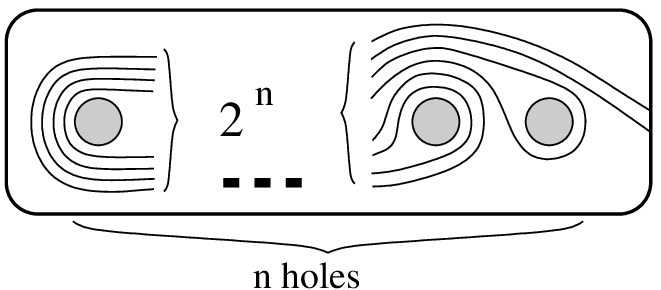,height=2.7cm}}
\begin{center}
                Figure 6.  
\end{center}

Denote the loops going around the
holes in $H_n'$ in the manner presented below by $g_1,g_2,...,g_n.$
\vspace*{.2in}\\

\centerline{\psfig{figure=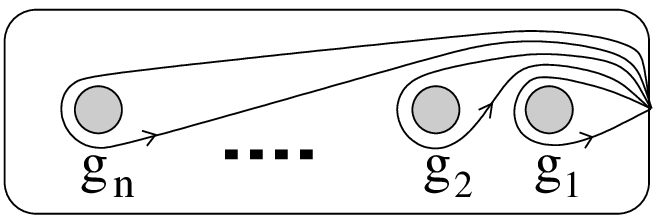,height=2cm}}
\begin{center}
                Figure 7. 
\end{center}

$F_n=\pi_1(H_n',x_0)$ is the free group on generators $g_1,...,g_n.$ 
We want to determine the presentation of the element $a_n\in F_n$
representing the arc $A$ given as in Fig. 6, with say an 
anti-clock orientation.
The presentations of $a_1$ and $a_2$ are as follows:
\vspace*{.1in}\\

\centerline{\psfig{figure=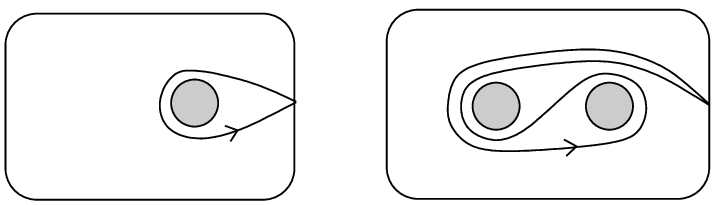,height=1.8cm}}
\begin{center}
                Figure 8: $a_1=g_1$ and $a_2=g_2g_1g_2^{-1}$ 
\end{center}

Observe that $a_{n+1}$ can be built inductively from $a_n$ by
replacing the $n$-th hole in $H_n'$ by two holes and twisting them
$180^0,$ see Fig 6. This operation corresponds to replacing 
all $g_n^{\pm 1}$ in the presentation of $a_n$ by
$g_{n+1}g_n^{\pm 1}g_{n+1}^{-1}.$ Therefore
$$a_3=g_3g_2g_3^{-1}g_1g_3g_2^{-1}g_3^{-1},\ a_4=g_4g_3g_4^{-1}g_2
g_4g_3^{-1}g_4^{-1}g_1g_4g_3g_4^{-1}g_2^{-1}g_4g_3^{-1}g_4^{-1}, 
{\rm \ etc.}$$

Note that each such presentation is reduced, i.e. none of the words $a_n$ 
has a subword of the form $g_i^{\pm 1}g_i^{\mp 1}.$ Indeed, this
is true for $n=1,2,3,4.$ In general, if this is true for some $n$ then it is
also true for $n+1:$ Suppose that $a_n$ is reduced. The only difference 
between $a_n$ and $a_{n+1}$ is that $g_n$ in $a_n$ is replaced by
$g_{n+1}g_ng_{n+1}^{-1}.$ Therefore, $a_{n+1}$ cannot contain
$g_i^{\pm 1}g_i^{\mp 1}$ for $i\leq n.$ Hence, if $a_{n+1}$ 
was not reduced, it would have to have a reduction in $g_{n+1}$'s. That is,
it would have to contain a subword composed of 
$g_{n+1}g_n^{\pm 1}g_{n+1}^{-1}$ and its inverse.
But this is impossible, since this would mean that $a_n$ contains the subword
$g_n^{\pm 1}g_n^{\mp 1}.$

Observe also that the number of appearances of $g_n^{\pm 1}$ in $a_n$
is $2^{n-1}.$ This can also be proved by induction: $g_1$ appears once
in $a_1,$ and $g_{n+1}$ appears in $a_{n+1}$ twice as many times as
$g_n$ in $a_n.$

Therefore we proved the following

\begin{proposition}\label{1.1}
The above inductively defined presentation of $a_n\in F_n$
is reduced and $g_n^{\pm 1}$ appears $2^{n-1}$ times in it.
\end{proposition}

Recall that we are supposed to count the minimal number of times the
the ring $S^1$ has to pass through the arc $A$ (in the original
puzzle) in a process of contracting it to a point.
Equivalently, we can calculate the number of times the arc
$A$ has to pass through the $n$-th hole in $H_n'$ in order to be
placed inside $H_{n-1}'\subset H_n'$ (in the situation in which the modified
puzzle is considered). We claim that this number is at least $2^{n-1}.$
 
Whenever $A$ passes through the $n$-th hole, $xg_n^{\pm
1}x^{-1}$ is inserted in a word representing $A$ or deleted from it.
Intuitively, since $g_n^{\pm 1}$ appears $2^{n-1}$ times in $a_n$
one has to repeat this process $2^{n-1}$ times. This may seem
as an obvious fact at first, but after a closer look one can realize that
it requires a proof.
We will see in the next section that this problem is a special case of 
an interesting problem in the combinatorial group theory.
In order to finish our proof we need the following fact.

\begin{proposition}\label{1.2}
If letters $b_1,b_2,...,b_k\in \{g_1^{\pm 1},...,g_{n-1}^{\pm 1}\}$ 
form a reduced word $b_1b_2...b_k$ then the minimal number of insertions 
or deletions of conjugates of $g_n^{\pm 1}$ into a word
$w= g_n^{\alpha_0}b_1g_n^{\alpha_1}b_2...b_kg_n^{\alpha_k}\in F_n$ 
necessary for transforming $w$ into an element of 
$F_{n-1}$ is $\sum_{i=0}^k |\alpha_i|$
\end{proposition}

We will prove the above proposition in the next section.
Right now we complete the proof of the main theorem:
Recall that $a_n$ can be obtained from $a_{n-1}$ by replacing all
$g_{n-1}^{\pm 1}$ appearing in $a_{n-1}$ by $g_ng_{n-1}^{\pm 1}g_n^{-1}.$
Therefore, $a_n= g_n^{\alpha_0}b_1g_n^{\alpha_1}b_2...b_kg_n^{\alpha_k},$
where, by Proposition \ref{1.1}, $b_1b_2...b_k=a_{n-1}$ is reduced and 
$\sum_{i=0}^k |\alpha_i|=2^{n-1}.$
By the last proposition one needs at least $2^{n-1}$ insertions of conjugates
of $g_n^{\pm 1}$ into a word representing $a_n$ or deletions of such 
conjugates from $a_n$ in order to eliminate all appearances of $g_n$ in it.
Therefore the arc $A$ has to pass at least $2^{n-1}$ times through the
$n$-th hole in $H_n'$ in order to be placed inside $H_{n-1}'\subset H_n'.$\\

\noindent {\bf Final Remark.} We found the complexity of the Chinese Rings
by the study of the homotopical properties of the arc $A$ in a
handlebody and by the use of the fundamental group. One can however consider
more sophisticated puzzles, whose solutions require more advanced tools
than homotopy theory. For example in order to show that the arc $A$ of
Fig. 9 can not be contracted to a point (without crossing changes) 
one can use the Kauffman bracket skein module of the solid torus [P] --
an algebraic construction which generalizes the Kauffman bracket polynomial,
[K].
\vspace*{.1in}\\

\centerline{\psfig{figure=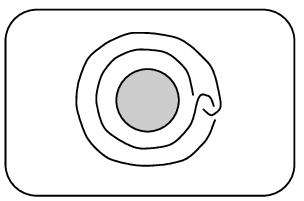,height=2.5cm}}

\begin{center}
                Figure 9 
\end{center}

\section{Problems in group theory involved in the proof of Kauffman's
conjecture}

The group theoretical problem which we encountered at the end of the
proof of the main theorem can be presented in a more general form
as follows: Consider a free group $F=<g_1,g_2,...\ |\ >$ and a quotient
map $\pi: F\to G = <g_1,g_2,...\ |\ r_1,r_2,...>.$
By an elementary operation on $w\in F$ we mean an insertion of a conjugate 
of $r^{\pm 1}_i$ into a word representing $w.$ Observe that
a deletion of $xr^{\pm 1}_ix^{-1}$ from a word representing
$w$ can be also realized as an insertion of $x^{-1}r^{\mp 1}_ix$ into this
word (followed by a reduction) and, therefore, it is also an
elementary operation. 
Observe also that two elements $x,y\in F$ have the same image, $\pi(x)=
\pi(y),$ if and only if they are related by a sequence of elementary
operations.

\noindent {\bf Problem 1}
Given $x,y\in F, \pi(x)=\pi(y),$ what is the minimal number of elementary
operations necessary for transforming $x$ into $y$? We denote this number by
$(x,y).$ If $\pi(x)\ne\pi(y)$ then we set $(x,y)=\infty.$

The above problem can be stated in a different form:
Any $w\in Ker(\pi)$ can be presented as $w= \prod_{i=1}^k 
p_ir_{\alpha_i}^{\pm 1}p_i^{-1}$, and we define $||w||$ to be 
the minimal number $k$ of conjugates of $r_{\alpha}^{\pm 1}$ 
appearing in such a presentation. If $w\not\in Ker(\pi)$ then $||w||=\infty.$

\noindent{\bf Problem 2}
Given $w\in Ker(\pi)$ calculate $||w||.$

It turns out that that Problems 1 and 2 are equivalent and that
the symbols $(\cdot,\cdot)$ and $||\cdot||$ have several interesting 
properties:

\begin{proposition}\label{2.1}
\begin{enumerate}
\item $||\cdot ||: F \to \{0,1,2,..,\infty\}$ is a ``norm'' on $F$ i.e.
(a) $||w||=0$ iff $w=e$ in $F,$ (b) $||vw||\leq ||v||+||w||,$
(c) $||w^{-1}||=||w||.$
\item $(x,y)=||xy^{-1}||$ and $(\cdot,\cdot)$ is a metric on $F,$
which is finite on $\pi^{-1}(g),$ for any $g\in G.$
\end{enumerate}
\end{proposition}

\begin{proof}
The proof of (1) is straightforward.
Also, the claim that $(\cdot,\cdot)$ is a metric on $F$ follows
immediately form (1) and the equality $(x,y)=||xy^{-1}||.$
Therefore we include only the proof of that equality:
If $\pi(x)\ne \pi(y)$ then $(x,y)=\infty=||xy^{-1}||.$ Hence we can
assume that $\pi(x)=\pi(y).$
Since $xy^{-1}=\prod_{i=1}^k p_ir_{\alpha_i}^{\pm 1}p_i^{-1},$ 
where $k=||xy^{-1}||,$ $x= \prod_{i=1}^k p_ir_{\alpha_i}^{\pm 1}p_i^{-1}\cdot
y$ can by transformed into $y$ by $k$ elementary operations. Therefore,
$||xy^{-1}||\geq (x,y).$

In order to prove the opposite inequality we need the following fact:\\
If $z$ is obtained from $1$ by $k$ elementary operations then $z$ can be
presented as $\prod_{i=1}^k p_ir_{\alpha_i}^{\pm 1}p_i^{-1}.$
We prove this fact by induction. The statement is true for $k=1.$
Suppose it is also true for $k-1$ and 
suppose that $z$ is obtained from $1$ by $k$ elementary operations.
By the inductive assumption, after $k-1$ operations we get 
a product of $k-1$ conjugates, $z'.$ The word $z$ is obtained by
inserting a word $qr_{\beta}^{\pm 1}q^{-1}$ into $z'.$ Therefore
$z=z_1qr_{\beta}^{\pm 1}q^{-1}z_2,$ for some $z_1,z_2$ such that $z'=z_1z_2.$
Thus $z=z'z_2^{-1}qr_{\beta}^{\pm 1}q^{-1}z_2$ is a product of $k$ conjugates.

Since $x$ is obtained from $y$ by $(x,y)$ elementary operations,
$xy^{-1}$ may be obtained from $yy^{-1}=1$ also by $(x,y)$ elementary
operations.
Therefore, by the fact proved above, $xy^{-1}$ can be presented as
a product of $(x,y)$ conjugates of $r_i^{\pm 1}$'s. Thus
$||xy^{-1}||\leq (x,y).$

\end{proof}

Notice that Problems 1-2 generalize the word problem in $G$ and that
they may be very difficult to solve in general. However,
fortunately to us, the group theoretic problem encountered at the end of 
the previous section is a special, solvable, case of Problems 1-2:\\
From now on $\pi:F_n\to F_{n-1}\subset F_n$ will be a projection given by 
$\pi(g_i)=g_i,$ for $i<n$ and $\pi(g_n)=1.$ Let $b_1b_2...b_k$ be a
reduced word in $F_{n-1}.$
We consider the problem of determining the minimal number of elementary
operations on $w=g_n^{\alpha_0}b_1g_n^{\alpha_1}b_2...b_kg_n^{\alpha_k}$
necessary for transforming this word into an element of $F_{n-1}.$
Since $w\in \pi^{-1}(b_1b_2...b_k),$ any word obtained by applying 
elementary operations to $w$ will still be in
$\pi^{-1}(b_1b_2...b_k).$ Observe that the only element of
$\pi^{-1}(b_1b_2...b_k)$ which can be presented as a word in the
letters $g_1^{\pm 1},...,g_{n-1}^{\pm 1}$
(i.e. a word without the letter $g_n^{\pm 1}$ in it) is
$b_1b_2...b_k\in F_{n-1}
\subset F_n.$ Therefore Proposition \ref{1.2} can be restated as\\

\noindent {\bf Proposition \ref{1.2}'}
{\em If letters $b_1,b_2,...,b_k\in\{g_1^{\pm 1},...,g_{n-1}^{\pm 1}\}$
form a reduced word $b_1b_2...b_k$ then}
$$(g_n^{\alpha_0}b_1g_n^{\alpha_1}b_2...b_kg_n^{\alpha_k},b_1b_2...b_k)=
\sum_{i=0}^k |\alpha_i|.$$

Let $P$ be a word of the form $s_1s_2...s_k$. Place points corresponding 
to $s_i$'s on the $x$-axis, $s_{i+1}$ after $s_i$; compare Fig. 10.
A connection, $C,$ on $P$ is a set of pairwise disjoint arcs in the
upper half plane. Each arc connects some $g_i$ with $g_{i}^{-1}$, 
in such a way that only some $g_n^{\pm 1}$'s may be left unconnected. 
The norm of the connection, $|C|$, is the number of unconnected letters. 
Fig. 10 shows two different connections for the word 
$P=g_1 g_n g_1^{-1}g_n g_1g_n^{-1}g_1^{-1}$, one of norm $3$
and another of norm $1.$\\

\centerline{\psfig{figure=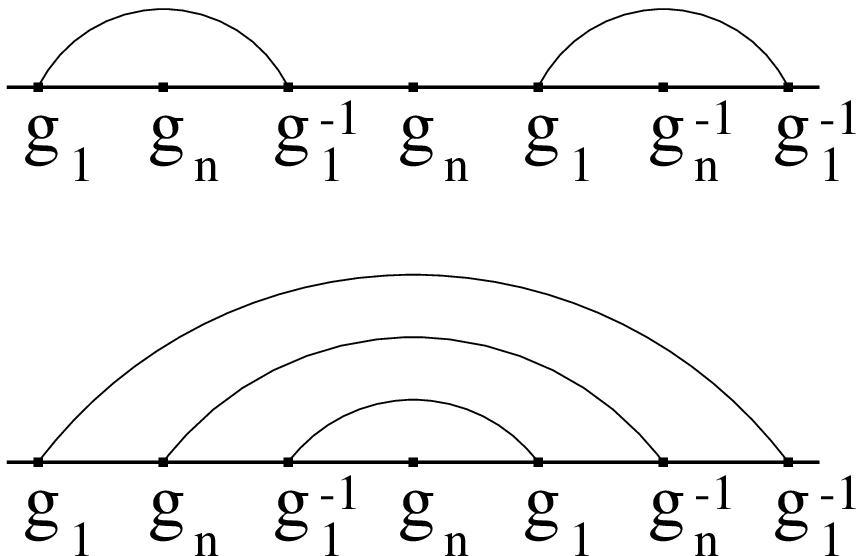,height=2.5cm}}

\begin{center}
                Figure 10 
\end{center}
We denote by $|w|_c$ the minimum of $|C|,$ where $C$ varies over
all connections on a word $w.$ If $w$ has no connection then
$|w|_c=\infty.$

\begin{theorem}\label{2.2}
For any word $w,$ $|w|_c=||w||.$
\end{theorem}

Observe, that since each word has only a finite number of connections
on it and it is easy to construct all of them,
the above theorem gives an explicit method of calculating $||w||.$

\begin{proof}
\begin{enumerate}
\item We prove $|w|_c\geq ||w||$ first.\\
If $w$ has no connection then $|w|_c=\infty$ and the inequality is
obvious. Therefore we may assume that $w$ has a connection $C.$
Consider two kinds of operations on $w$:\\
(I) If there is an unconnected letter $g_n^{\pm 1}$ in $w$ then we 
delete this letter and we obtain a new word $w'$ with a connection $C'$
(composed of the same arcs as $C$).\\ 
(II) If there is a pair $g_i,g_i^{-1}$ of letters in $w$ connected
by an arc which is not nested (i.e. $g_i,g_i^{-1}$ are neighbors) then
we remove these letters and the arc connecting them and we obtain a new
word $w'$ of a shorter length with a connection $C'.$\\
Observe that each of the above operations decreases the length of $w,$ 
and we can always apply at least one of them to $w,$ unless $w$ is the
trivial word $1.$ Therefore any word $w$ can be reduced to $1$ by a 
sequence of operations of the first and the second type.
Observe that Operation I changes $||w||$ by at most one and decreases 
the norm of the connection on $w$ by $1.$ Operation II does not change 
$||w||$ nor the norm of the connection on $w.$ Since at the end of the process
(when $w=1$) both $||w||$ and the norm of the connection on $w$ are $0,$
$|C|\geq ||w||.$ Since this inequality holds for any connection $C$ 
on $w,$ we also have $|w|_c\geq ||w||.$

\item We also claim that $|w|_c\leq ||w||.$\\
We can assume that $w\in Ker(\pi),$ since otherwise $||w||=\infty.$
There is a word $w'= \prod_{i=1}^{k} 
p_ig_n^{\pm 1}p_i^{-1}$ representing the same element of the
group $F_n$ as the word $w,$ with $k=||w||$. Each factor 
$p_ig_n^{\pm 1}p_i^{-1}$ has a connection of a norm $1$ (a nested family
of arcs). Therefore the product $w'$ has a connection of a norm $||w||$. 
The word $w'$ can be transformed into the word $w$ 
by a sequence of insertions and deletions of subwords of the form
$g_i^{\pm 1}g_i^{\mp 1},\ i=1,2,...,n.$
Observe that after each insertion we obtain a new word with a
connection of the same norm. Moreover, it is not difficult to see, that
after each deletion we obtain a new
word with a connection of the same norm, if $i<n,$ or lower or equal
norm, if $i=n.$
Therefore $w$ has a connection $C$ of norm $|C|\leq ||w||$ and, hence,
$|w|_c\leq ||w||.$
\end{enumerate}
\end{proof}

{\it Proof of Proposition \ref{1.2}':}
Since 
$$(g_n^{\alpha_0}b_1g_n^{\alpha_1}b_2...b_kg_n^{\alpha_k},b_1b_2...b_k)\leq
\sum_{i=0}^k |\alpha_i|,$$
we only need to prove the opposite inequality.
By Proposition \ref{2.1} and Theorem \ref{2.2} it is enough to prove
that each connection on
$$g_n^{\alpha_0}b_1g_n^{\alpha_1}b_2...b_kg_n^{\alpha_k}b_k^{-1}...b_2^{-1}
b_1^{-1}$$ has its norm greater or equal to $\sum_{i=0}^k |\alpha_i|.$
Let $C$ be any such connection. Observe that it is enough to prove
that $C$ does not connect any pair of letters $g_n^{\pm 1}, g_n^{\mp 1}.$
Suppose that $C$ does connect letters
$g_n^{\pm 1}$ and $g_n^{\mp 1}$ enclosing a word
$b_ig_n^{\alpha_i}b_{i+1}...g_n^{\alpha_{j-1}}b_j.$ The connection $C$ 
restricts to a connection on this word and therefore 
$b_ig_n^{\alpha_i}b_{i+1}... g_n^{\alpha_{j-1}}b_j\in Ker(\pi).$ Hence 
$b_ib_{i+1}...b_j=e$ in $F_{n-1}=<g_1,...,g_{n-1}>$ and therefore 
$b_1b_2...b_k$ is reducible, what contradicts our assumption.
\Box\\

\vspace{.2in}\ \\

\noindent Authors address:\\
\centerline{\it Department of Mathematics, University of Maryland}
\centerline{\it College Park, MD 20742}
\centerline{\it e-mails: przytyck@gwu.edu and asikora@math.umd.edu}
\vspace*{.2cm}

\noindent The first author is on leave from\\
\centerline{\it Department of Mathematics, The George Washington University}
\centerline{\it Washington, DC 20052}

\begin{thebibliography}{99}
\bibitem [BC] {BC} W.W. Rouse Ball and H.S.M. Coxeter, 
Mathematical Recreations and Essays, University of Toronto Press, 
Toronto, 1974 (First edition - 1892).
\bibitem [BL] {BL} J. S. Birman, X. L. Lin, Knot polynomials and Vassiliev's
invariants, {\em Invent. Math.}, {\bf 111} (1993), 225--270.
\bibitem [H] {Ha} K. Habiro, Claspers and finite type invariants of links,
{\it Geom. Topol.}\\ {\bf 4} (2000), 1--83,
http://www.maths.warwick.ac.uk/gt/GTVol4/\allowbreak paper1.abs.html
\bibitem [K1] {K1} L. H. Kauffman, An invariant of regular isotopy, 
{\em Trans. Amer. Math. Soc.} {\bf 318} (1990), no. 2, 417--471.
\bibitem [K2] {K2}L. H. Kauffman, Tangle complexity and
the topology of the Chinese rings, in Mathematical approaches
to biomolecular structure and dynamics, IMA
Vol. 82, Springer, New York, 1996, 1--10.
\bibitem [P] {P} J. H. Przytycki, Fundamentals of Kauffman bracket skein
modules, {\it Kobe J. Math.}, 16(1), 1999, 45-66,
http://xxx.lanl.gov/abs/math.\allowbreak GT/9809113
\end{thebibliography}
\end{document}